\documentclass[reqno]{amsart}
\usepackage[english]{babel}
\usepackage{amssymb,verbatim}
\usepackage[T1,T5]{fontenc}

\theoremstyle{definition}

\theoremstyle{remark}

\numberwithin{equation}{section}

\newcommand{\cH}{\mbox{$\mathcal H$}}

\newcommand{\n}{\mbox{\noindent}}
\def\interior{\mathop{\rm int}\nolimits}

\begin{document}

\author{Pham Hoang Hiep}
\address{Institute of Mathematics\\
Vietnam Academy of Science and Technology \\ 
18, Hoang Quoc Viet, Hanoi, Vietnam} 
\email{phhiep@math.ac.vn}

\title{Log canonical thresholds and Monge-Amp\`ere masses}

\begin{abstract}
In this paper, we prove an inequality for log canonical thresholds and Monge-Amp\`ere masses. The idea of proof is a combination of the Ohsawa-Takegoshi $L^2$-extension theorem and inequalities in \cite{ACKHZ09}, \cite{DH14}. It also give an analytic proof for the main result in \cite{DH14} in the case of dimension 2.
\end{abstract}

\maketitle

\section{Notation and main results}

\n

Here we put $d^c=\frac {i} {2\pi} (\overline {\partial} - \partial)$, so that $dd^c=\frac {i} {\pi}\partial \overline {\partial}$. The normalization of the $d^c$ operator is chosen so that we have precisely $(dd^c\log|z|)^n=\delta_0$ for the Monge-Amp\`ere operator in $\mathbb C^n$. The Monge-Amp\`ere operator is defined on locally bounded plurisubharmonic functions according to the definition of Bedford-Taylor (\cite{BT76}, \cite{BT82}); it can also be extended to plurisubharmonic functions with isolated or compactly supported poles by \cite{Dem93}.

If $\Omega$ is an open subset of $\mathbb C^n$, we let $PSH(\Omega)$ (resp. $PSH^-(\Omega)$) denote the set of plurisubharmonic (resp. negative plurisubharmonic) functions on $\Omega$.

\medskip

\n{\bf Definition 1.1.} Let $\Omega$ be a bounded hyperconvex domain (i.e. a domain possessing a negative plurisubharmonic exhaustion function). Following Cegrell (\cite{Ceg98}, \cite{Ceg04}), we introduce certain classes of plurisubharmonic functions on $\Omega$, in relation with the definition of the Monge-Amp\`ere operator $:$

$$\mathcal E_0(\Omega)=\{\varphi\in PSH^-(\Omega):\ \lim\limits_{z\to\partial\Omega}\varphi (z)=0,\ \int\limits_\Omega (dd^c\varphi)^n<+\infty\},$$

$$\mathcal E_p(\Omega)=\{\varphi\in PSH^-(\Omega):\ \exists\ \mathcal E_0(\Omega)\ni\varphi_j\searrow\varphi,\ \sup\limits_{j\geq 1}\int\limits_\Omega (-\varphi_j )^p (dd^c\varphi_j)^n<+\infty\},$$

$$\mathcal F(\Omega)=\{\varphi\in PSH^-(\Omega):\ \exists\ \mathcal E_0(\Omega)\ni\varphi_j\searrow\varphi,\ \sup\limits_{j\geq 1}\int\limits_\Omega(dd^c\varphi_j)^n<+\infty\},$$

$$\mathcal E(\Omega)=\{\varphi\in PSH^-(\Omega):\ \exists\ \varphi_K\in\mathcal F(\Omega)\ \text{such that}\ \varphi_K=\varphi\ \text{on}\ K,\ \forall K\subset\subset\Omega\}.$$

It is proved in \cite{Ceg04} that the class $\mathcal E(\Omega )$ is the biggest subset of $PSH^-(\Omega)$ on which the Monge-Amp\`ere operator is well-defined (see also \cite{Blo06}). It is easy to see that $\smash{\mathcal E(\Omega )}$ contains the class of negative plurisubharmonic functions which are locally bounded outside isolated singularities.

For $\varphi\in PSH(\Omega)$ and $0\in\Omega$, we introduce the log canonical threshold at $0$

$$c(\varphi) = \sup\big\{ c>0:\ e^{-2c\varphi} \text{ is } L^1 \text{ on a neighborhood of }0\big\}.$$
It is an invariant of the singularity of $\varphi$ at $0$. We refer to \cite{Ber15}, \cite{DK01}, \cite{GZ15a}, \cite{GZ15b}, \cite{GZ15c}, \cite{Hie14}, \cite{Kim15}, \cite{Kis87}, \cite{Mat15}, \cite{Kis94} and \cite{Ras14} for further information and applications to this number.

For $1\leq k\leq n$, we denote
$$c_k (\varphi) = \sup \{ c( \varphi|H ): \text{ H is a subspace of k-dimension through }0\},$$
where $\varphi|H$ is the restriction of $\varphi$ to the subspace $H$. For $\varphi\in\mathcal E(\Omega)$ and $1\leq k\leq n$, we introduce the intersection numbers
$$e_{k} (\varphi ) = \int_{\{0\}} (dd^c\varphi)^k\wedge (dd^c\log\Vert z\Vert)^{n-k}$$
which can be seen also as the Lelong numbers of $(dd^c\varphi)^k$ at $0$.

Our main result is the following estimate. It is a generalization of inequalities in \cite{ACKHZ09}, \cite{Cor95}, \cite{Cor00}, \cite{DH14}, \cite{FEM03}, \cite{FEM04}; such inequalities have fundamental applications to birational geometry (see \cite{IM72}, \cite{Puk87}, \cite{Puk02}, \cite{Isk01}, \cite{Che05}). 

\medskip

{\bf Theorem 1.2.} {\it Let $\varphi\in\mathcal E(\Omega)$ and $0\in\Omega$. Then $c(\varphi)=+\infty$

if $e_1(\varphi)=0$, and otherwise

$$c(\varphi)\geq c_{n-1} (\varphi ) + \frac { (n-1)^{n-1} } { c_{n-1} (\varphi )^{n-1} e_n(\varphi ) }.$$}

\medskip

{\bf Remark 1.3.} i) Theorem 1.2 also give an analytic proof for the main result in \cite{DH14} in the case of dimension 2. 

ii) In \cite{GZ15c}, Quan Guan and Xiangyu Zhou showed that
$$c_k (\varphi) - c_{k-1} (\varphi)\leq c_{k-1} (\varphi) - c_{k-2} (\varphi).$$
By induction on the dimension $n$, we obtain an upper bound for the log canonical thresholds
$$c(\varphi )\leq \frac {n} {n-1} c_{n-1} (\varphi ).$$

\section{Some lemmas}

\n

We start with a monotonicity statement.
\medskip

\n

{\bf Lemma 2.1.} {\it Let $\varphi,\psi\in\mathcal E(\Omega )$ be such that $\varphi \leq \psi$ $($i.e.\ $\varphi$ is ``more singular'' than $\psi)$. Then

$$c_{n-1} (\varphi ) + \frac { (n-1)^{n-1} } { c_{n-1} (\varphi )^{n-1} e_n(\varphi ) } \leq c_{n-1} (\psi ) + \frac { (n-1)^{n-1} } { c_{n-1} (\psi )^{n-1} e_n(\psi ) }.$$}

\n {\it Proof.} First, by the main inequality of \cite{FEM03}, \cite{FEM04}, \cite{Dem09} and Corollary 2.2 in \cite{DH14}, we have

$$c_{n-1} (\varphi ) \geq \frac { n-1 } { e_{n-1}(\varphi )^{ \frac {1} {n-1} } }\geq \frac { n-1 } { e_{n}(\varphi )^{ \frac {1} {n} } }.$$

Now, we set

$$D=\{t=(t_1,t_2)\in [0,+\infty)^2:\ (n-1)^{n}t_1^{n} \leq t_2\}.$$

Then $D$ is a convex set in $[0,+\infty)^2$ and $(\frac {1} {c_{n-1} (\varphi ) } ,e_n (\varphi)), (\frac {1} {c_{n-1} (\psi ) } ,e_n (\psi ))\in D$. We consider the function $f: \interior D\to [0,+\infty)$.

$$f(t_1,t_2) = \frac {1} {t_1} + \frac {(n-1)^{n-1} t_1^{n-1} } { t_2 }.$$

We have

$$\frac {\partial f} {\partial t_1} (t) = -\frac { 1 } { t_1^2 } + \frac {(n-1)^{n} t_1^{n-2} } { t_2 }\leq 0,\ \forall t\in D;$$

$$\frac {\partial f} {\partial t_2} (t) = -\frac {(n-1)^{n-1} t_1^{n-1} } { t_2^2 }\leq 0,\ \forall t\in D.$$

For $a,b\in\interior D$ such that $a_j\geq b_j$, $\forall j=1,2$, $[0,1]\ni \lambda\to f( b+ \lambda(a-b) )$ is a decreasing function. This implies that $f(a)\leq f(b)$ for $a,b\in\interior D$, $a_j\geq b_j$, $\forall j=1,2$. On the other hand, the hypothesis $\varphi\leq\psi$ implies $\frac {1} { c_{n-1} (\varphi ) } \geq \frac {1} { c_{n-1} (\psi ) }$ and $e_n (\varphi) \geq e_n (\psi)$, by the comparison principle \cite{Dem87}. Therefore $f(\frac {1} {c_{n-1} (\varphi ) } ,e_n (\varphi))\leq f(\frac {1} {c_{n-1} (\psi ) }  ,e_n (\psi))$.\qed

\medskip

Next, we need a following lemma which follows from \cite{GZ15a} and \cite{Hie14}.

\medskip

{\bf Lemma 2.2.} {\it Let $\varphi \in PSH^- (\Delta^n )$ and $c>0$ be such that
$$\varliminf\limits_{ z_n\to 0 }\int\limits_{\Delta^{n-1} } e^{-2c\varphi (z',z_n )} dV_{n-2} (z') |z_n|^2 = 0,$$
where $\Delta$ is the unit disk in $\mathbb C$. Then
$$c(\varphi )\geq c.$$}

\n {\it Proof.}
We set
$$h(z_n) = \int_{\Delta^{n-1} } e^{-2c\varphi (z',z_n )} dV_{n-2} (z')$$
By the $L^2$-extension theorem of Ohsawa and Takegoshi (\cite{OT87}, see also \cite{Blo13} and \cite{Dem15}), there exists a holomorphic function $f$ on $\Delta^{n}$ such that $f(z',w_n)=1$ for all $z'\in\Delta^{n-1}$,
and
$$\int\limits_{\Delta^n } |f(z)|^2 e^{-2c\varphi (z)} d V_{2n}(z)\leq \frac {1} {\pi} \int\limits_{\Delta^{n-1}} e^{-2 c\varphi (z',w_n) } d V_{2n-2}(z') = \frac { h(w_n) } { \pi },$$
By the mean value inequality for the plurisubharmonic function $|f|^2$, we get
\begin{eqnarray*}
|f(z)|^2&\leq& \frac 1 {\pi^n(1 - |z_1|)^2\ldots(1 - |z_n|)^2}\int\limits_{ \Delta_{1 - |z_1|} (z_1)\times\ldots\times\Delta_{1 - |z_n|} (z_n) } |f|^2 dV_{2n}\\
&\leq& \frac {h(w_n)}{ \pi^{n+1}(1 - |z_1|)^2\ldots(1 - |z_n|)^2 },
\end{eqnarray*}
where $\Delta_r(z)$ is the disc of center $z$ and radius $r$. Hence, for any $r<\frac 1 2$, we have
\begin{equation}
\Vert f\Vert _{ L^\infty (\Delta_{r}^n) }\leq \frac { h(w_n)^{\frac 1 2 } }{\pi^{\frac {n+1} 2} (1-r)^n }.
\end{equation}
Since $f(z',w_n)-1=0$, $\forall z'\in\Delta^{n-1}$, we can write $f(z)=1+(z_n-w_n)g(z)$ for some function holomorphic function $g$ on $\Delta^n$. By (2.1), for any $|w_n|<r<\frac 1 2$, we get
\begin{eqnarray*}
\Vert g\Vert _{ \Delta_{r}^n } = \Vert g\Vert _{ \Delta_{r}^{n-1}\times\partial\Delta_{r} }
&\leq& \frac {1} {r-|w_n|} \Big(\Vert f\Vert _{ L^\infty (\Delta_{r}^n) }+1\Big)\\
&\leq& \frac {1} {r-|w_n|} \Big(\frac { h(w_n)^{\frac 1 2 } }{\pi^{\frac {n+1} 2} (1-r)^n }+1\Big).
\end{eqnarray*}
Hence
$$|f(z)-1| = |z_n - w_n||g(z)|\leq |z_n - w_n| \Vert g\Vert _{ \Delta_{r}^n } \leq \frac {2|w_n|} {r-|w_n|} \Big(\frac { h(w_n)^{\frac 1 2 } }{\pi^{\frac {n+1} 2} (1-r)^n }+1\Big),$$
for all $(z',z_n)\in\Delta_r^{n-1}\times\Delta_{|w_n|}$. Take $w_n\in\Delta_r$ such that
$$|f(z)-1| \leq \frac {2|w_n|} {r-|w_n|} \Big(\frac { h(w_n)^{\frac 1 2 } }{\pi^{\frac {n+1} 2} (1-r)^n }+1\Big) < 1,$$
for all $(z',z_n)\in\Delta_r^{n-1}\times\Delta_{|w_n|}$. Combination these above inequalities, we get that
$$\int\limits_{\Delta_r^{n-1}\times\Delta_{|w_n|} } e^{-2c\varphi (z)} d V_{2n}(z) < +\infty.$$
This implies that $c(\varphi )\geq c$.  \qed

Next, we need a lemma which follows from \cite{ACKHZ09}.

{\bf Lemma 2.3.} {\it Let $\varphi\in\mathcal E_1 (\Omega )$. Then
\begin{equation*}
V_{2 n} (\{ \varphi < -t \}) \leq c_n \delta _\Omega^{2 n} \left(1 + t^{ \frac {n+1} {n} } ( \int\limits_{\Omega} -\varphi (dd^c\varphi )^n )^{- \frac 1 n} \right)^{n - 1} e^ { -  2 n  t^{ \frac {n+1} {n} } ( \int\limits_{\Omega} -\varphi (dd^c\varphi )^n )^{- \frac 1 n} },
\end{equation*}
where  $\delta _\Omega := diam (\Omega)$ is the diameter of $\Omega$ and $c_n$ is a constant.}

\n {\it Proof.} Using Proposition 6.1 in \cite{ACKHZ09}, we have
\begin{equation*}
V_{2 n} (\{ \varphi < -t \}) \leq c_n \delta _\Omega^{2 n} \left(1 + Cap_{\Omega} ( \{ \varphi < -t \} )^{- \frac 1 n} \right)^{n - 1} e^ {-  2 n \ Cap_{\Omega} (\{ \varphi < -t \})^{- \frac 1 n} },
\end{equation*}
On the other hand, since $t^{-1} \varphi\leq h_{\{ \varphi < -t \},\Omega }$, we have
$$Cap_{\Omega} ( \{ \varphi < -t \} ) = \int\limits_\Omega -h_{\{ \varphi < -t \},\Omega } (dd^c h_{\{ \varphi < -t \},\Omega } )^n\leq t^{-n-1}\int\limits_{\Omega} -\varphi (dd^c\varphi )^n.$$
Therefore
\begin{equation*}
V_{2 n} (\{ \varphi < -t \}) \leq c_n \delta _\Omega^{2 n} \left(1 + t^{ \frac {n+1} {n} } ( \int\limits_{\Omega} -\varphi (dd^c\varphi )^n )^{- \frac 1 n} \right)^{n - 1} e^ { -  2 n  t^{ \frac {n+1} {n} } ( \int\limits_{\Omega} -\varphi (dd^c\varphi )^n )^{- \frac 1 n} }.
\end{equation*}
\qed

\n
{\bf Lemma 2.4.} {\it Let $\varphi\in\mathcal E_1 (\Omega )$ and $\Omega'\subset\Omega$ be such that
$$\int\limits_{\Omega } -\varphi (dd^c \varphi )^{n} \leq A,$$
and
$$\int\limits_{\Omega' } e^{ -2c\varphi } dV_{2n} \leq B.$$
Then for $\lambda\in (c,c+\frac c n)$
$$\aligned\int\limits_{\Omega' } e^{ -2\lambda \varphi } dV_{2n}&\leq V_{2n} (\Omega') + \frac { B } { 2(\lambda-c) } e^{ 2(\lambda - c) A c^{n} n^{-n}  } \\
&\ + \frac {1} {2} c_n \delta _\Omega^{2 n} \left(1 + 2^{n+1} A \lambda^{n+1} n^{-n-1} \right)^{n - 1} ( (n+1) c n^{-1} -\lambda )^{-1} e^{ 2(\lambda - c) A c^{n} n^{-n}  }\\
&\ + c_n \delta _\Omega^{2 n} \frac { n^2 } { \lambda (n+2) } \int\limits_{ 0 }^{ +\infty }  \left( 1 + 2x  \right)^{n - 1} e^ { -  2 n  x } dx,
\endaligned$$
where  $\delta _\Omega := diam (\Omega)$ is the diameter of $\Omega$ and $c_n$ is a constant.}

\n {\it Proof.} By Lemma 2.3, we have
\begin{equation*}
V_{2 n} (\{ \varphi < -t \}) \leq c_n \delta _\Omega^{2 n} \left(1 + t^{ \frac {n+1} {n} } A^{- \frac 1 n} \right)^{n - 1} e^ { -  2 n  t^{ \frac {n+1} {n} } A^{- \frac 1 n} }.
\end{equation*}
On the other hand, we have
$$V_{2n} (\{ \varphi < -t \}\cap\Omega') \leq e^{ -2 c t } \int\limits_{\Omega' } e^{ -2c\varphi } dV_{2n} \leq B e^{ -2 c t }.$$
Therefore
$$V_{2n} (\{ \varphi < -t \}\cap\Omega')\leq \min ( c_n \delta _\Omega^{2 n} \left(1 + t^{ \frac {n+1} {n} } A^{- \frac 1 n} \right)^{n - 1} e^ { -  2 n  t^{ \frac {n+1} {n} } A^{- \frac 1 n} }, B e^{ -2 c t } ).$$
This implies that

\n
$\int\limits_{\Omega' } e^{ -2\lambda\varphi } dV_{2n} = V_{2n} (\Omega') + \int\limits_0^{+\infty } e^{2\lambda t} V_{2 n} (\{ \varphi < -t \}\cap\Omega' ) dt$

\n
$\leq V_{2n} (\Omega') + B\int\limits_0^{ Ac^{n} n^{-n}  } e^{2(\lambda - c) t} dt$

\n
$+ \int\limits_{ Ac^{n} n^{-n}  }^{+\infty } c_n \delta _\Omega^{2 n} \left(1 + t^{ \frac {n+1} {n} } A^{- \frac 1 n} \right)^{n - 1} e^ { -  2 n  t^{ \frac {n+1} {n} } A^{- \frac 1 n} + 2\lambda t} dt$

\n
$\leq V_{2n} (\Omega') + \frac { B } { 2(\lambda-c) } e^{ 2(\lambda - c) A c^{n} n^{-n}  }$

\n
$+ c_n \delta _\Omega^{2 n} \int\limits_{ Ac^{n} n^{-n}  }^{+\infty }  \left(1 + t^{ \frac {n+1} {n} } A^{- \frac 1 n} \right)^{n - 1} e^ { -  2 n  t^{ \frac {n+1} {n} } A^{- \frac 1 n} + 2\lambda t} dt$

\n
$\leq V_{2n} (\Omega') + \frac { B } { 2(\lambda-c) } e^{ 2(\lambda - c) A c^{n} n^{-n}  }$

\n
$+ c_n \delta _\Omega^{2 n} \int\limits_{ Ac^{n} n^{-n}  }^{ 2^n A \lambda^n n^{-n} }  \left(1 + t^{ \frac {n+1} {n} } A^{- \frac 1 n} \right)^{n - 1} e^ { -  2 n  t^{ \frac {n+1} {n} } A^{- \frac 1 n} + 2\lambda t} dt$

\n
$+c_n \delta _\Omega^{2 n} \int\limits_{ 2^n A \lambda^n n^{-n}  }^{ +\infty }  \left(1 + t^{ \frac {n+1} {n} } A^{- \frac 1 n} \right)^{n - 1} e^ { -  2 n  t^{ \frac {n+1} {n} } A^{- \frac 1 n} + 2\lambda t} dt$

\n
$\leq V_{2n} (\Omega') + \frac { B } { 2(\lambda-c) } e^{ 2(\lambda - c) A c^{n} n^{-n}  } + I_1 + I_2,$

\n
where
$$I_1 = c_n \delta _\Omega^{2 n} \left(1 + 2^{n+1} A \lambda^{n+1} n^{-n-1} \right)^{n - 1} \int\limits_{ Ac^{n} n^{-n}  }^{ 2^n A \lambda^n n^{-n} }  e^ { -  2 n ( t^{ \frac {n+1} {n} } A^{- \frac 1 n} - \lambda n^{-1} t )} dt;$$
and
$$I_2 =  c_n \delta _\Omega^{2 n} \int\limits_{ 2^n A \lambda^n n^{-n}  }^{ +\infty }  \left(1 + 2 ( t^{ \frac {n+1} {n} } A^{- \frac 1 n} -\lambda n^{-1} t ) \right)^{n - 1} e^ { -  2 n  ( t^{ \frac {n+1} {n} } A^{- \frac 1 n} - \lambda n^{-1} t) } dt.$$

\n
Set $x = t^{ \frac {n+1} {n} } A^{- \frac 1 n} -\lambda n^{-1} t$. We have
$$x'(t) = n^{-1} ( (n+1) t^{\frac 1 n} A^{-\frac 1 n} -\lambda )\geq n^{-1} ( (n+1) c n^{-1} -\lambda ) > 0,\ \forall t > Ac^{n} n^{-n},$$
and
$$x'(t)\geq \lambda (n+2) n^{-2}\ \forall t > 2^n A \lambda^n n^{-n}.$$
Using the change of variable $x = t^{ \frac {n+1} {n} } A^{- \frac 1 n} -\lambda n^{-1} t$, we get
$$\aligned I_1 &= c_n \delta _\Omega^{2 n} \left(1 + 2^{n+1} A \lambda^{n+1} n^{-n-1} \right)^{n - 1} \int\limits_{ Ac^{n} n^{-n}  }^{ 2^n A \lambda^n n^{-n} }  e^ { -  2 n ( t^{ \frac {n+1} {n} } A^{- \frac 1 n} - \lambda n^{-1} t )} dt\\
&=c_n \delta _\Omega^{2 n} \left(1 + 2^{n+1} A \lambda^{n+1} n^{-n-1} \right)^{n - 1} n ( (n+1) c n^{-1} -\lambda )^{-1} \int\limits_{ Ac^{n} n^{-n-1} (c-\lambda )  }^{ 2^n A \lambda^{n+1} n^{-n-1} } e^{-2n x} dx\\
&= \frac {1} {2} c_n \delta _\Omega^{2 n} \left(1 + 2^{n+1} A \lambda^{n+1} n^{-n-1} \right)^{n - 1} ( (n+1) c n^{-1} -\lambda )^{-1} e^{ 2(\lambda - c) A c^{n} n^{-n}  };
\endaligned$$
and
$$\aligned I_2 &= c_n \delta _\Omega^{2 n} \int\limits_{ 2^n A \lambda^n n^{-n}  }^{ +\infty }  \left(1 + 2 ( t^{ \frac {n+1} {n} } A^{- \frac 1 n} -\lambda n^{-1} t ) \right)^{n - 1} e^ { -  2 n  ( t^{ \frac {n+1} {n} } A^{- \frac 1 n} - \lambda n^{-1} t) } dt\\
&=c_n \delta _\Omega^{2 n} \frac { n^2 } { \lambda (n+2) } \int\limits_{ 2^n A \lambda^{n+1} n^{-n-1}  }^{ +\infty }  \left( 1 + 2x  \right)^{n - 1} e^ { -  2 n  x } dx\\
&= c_n \delta _\Omega^{2 n} \frac { n^2 } { \lambda (n+2) } \int\limits_{ 0 }^{ +\infty }  \left( 1 + 2x  \right)^{n - 1} e^ { -  2 n  x } dx.
\endaligned$$
Combining above inequalities, we get that
$$\aligned\int\limits_{\Omega' } e^{ -2\lambda \varphi } dV_{2n}&\leq V_{2n} (\Omega') + \frac { B } { 2(\lambda-c) } e^{ 2(\lambda - c) A c^{n} n^{-n}  } \\
&\ + \frac {1} {2} c_n \delta _\Omega^{2 n} \left(1 + 2^{n+1} A \lambda^{n+1} n^{-n-1} \right)^{n - 1} ( (n+1) c n^{-1} -\lambda )^{-1} e^{ 2(\lambda - c) A c^{n} n^{-n}  }\\
&\ + c_n \delta _\Omega^{2 n} \frac { n^2 } { \lambda (n+2) } \int\limits_{ 0 }^{ +\infty }  \left( 1 + 2x  \right)^{n - 1} e^ { -  2 n  x } dx.
\endaligned$$  \qed

\n

\section{Proof of the main theorem}

{\bf 3.1. Reduction to the case of plurisubharmonic functions with analytic singularities.}

\n

In the first step, we reduce the proof to the case $\varphi = \log (|f_1|^2+\ldots+|f_N|^2)$, where $f_1,\ldots,f_N$ are germs of holomorphic functions at $0$. Following \cite{Dem92}, we let $\cH_{m\varphi}(\Omega)$ be the Hilbert space of holomorphic functions $f$ on $\Omega$ such that
$$\int_\Omega|f|^2e^{-2m\varphi}dV<+\infty,$$
and let $\psi_m=\frac{1}{2m}\log\sum|g_{m,k}|^2$ where $\{g_{m,k}\}_{k\geq 1}$ is an orthonormal basis of $\cH_{m\varphi}(\Omega)$. Thanks to Theorem 4.2 in \cite{DK01}, mainly based on the Ohsawa-Takegoshi $L^2$ extension theorem \cite{OT87} (see also \cite{Dem92}), there are constants $C_1,C_2>0$ independent of $m$ such that

$$\varphi(z)-\frac{C_1}{m}\leq \psi_m(z) \leq \sup_{ |\zeta-z| < r}\varphi(\zeta )+\frac{1}{m}\log\frac{C_2}{r^n},$$
for every $z\in\Omega$ and $r<d(z,\partial\Omega)$. Moreover
$$\nu(\varphi)-\frac{n}{m}\leq\nu(\psi_m)\leq\nu(\varphi),$$

$$\frac 1 {c(\varphi)}-\frac{1}{m}\leq\frac{1}{c(\psi_m)}\leq \frac{1}{c(\varphi)}.$$

By Lemma 2.1, we have

$$f(\frac {1} {c_{n-1} (\varphi ) } ,e_n (\varphi))\leq f(\frac {1} {c_{n-1} (\psi_m ) } ,e_n (\psi_m )),\qquad \forall m\geq 1.$$

The above inequalities show that in order to prove the lower bound of $c(\varphi)$ in Theorem 1.2, we only need to prove it for $c(\psi_m)$ and let $m$ tend to infinity. Also notice that since the Lelong numbers of a function $\varphi\in{\mathcal E}(\Omega)$ occur only on a discrete set, the same is true for the functions $\psi_m$.

\n
{\bf 3.2. Reduction to the case $\varphi\in\mathcal F(\Delta^n ) \cap C(\Delta^n\backslash \{ 0 \})$ with $c_{n-1} (\varphi ) = c(\varphi (z',0) )$ and $(dd^c\varphi )^n = \delta_{ \{0\} }$, where $\delta_{ \{0\} }$ is the Dirac measure at the point $\{0\}$.}

In the second step, we reduce the proof to the case $\varphi\in\mathcal F(\Delta^n ) \cap C(\Delta^n\backslash \{ 0 \})$ and $c_{n-1} (\varphi ) = c(\varphi (z',0) )$. By first step, we can assume that $\varphi\in PSH^-(\Delta^n ) \cap C^\infty (\Delta^n\backslash \{ 0 \} )$. For $0<r<1$, we set
$$\psi_r = \sup\{ u\in PSH^-(\Delta^n ):\ u\leq\varphi \text{ on } \Delta_r^n \}.$$
Then $\psi_r\in\mathcal F(\Delta^n ) \cap C(\Delta^n\backslash \{ 0 \})$ and $\psi_r\nearrow \psi\in \mathcal F(\Delta^n ) \cap C(\Delta^n\backslash \{ 0 \})$ as $r\searrow 0$. We have $(dd^c \psi )^n = \int\limits_{\Omega} (dd^c \psi )^n \delta_{ \{0\} }$ and $c(\varphi ) = c(\psi )$. By Lemma 2.1, we have

$$f(\frac {1} {c_{n-1} (\varphi ) } , e_n (\varphi))\leq f(\frac {1} {c_{n-1} (\psi ) } , e_n (\psi )).$$

The above inequalities show that in order to prove the lower bound of $c(\varphi)$ in Theorem 1.2, we only need to prove it for $c(\psi)$. 

{\bf 3.3. Proof of the main theorem in the case $\varphi\in\mathcal F(\Delta^n ) \cap C(\Delta^n\backslash \{ 0 \})$ with $c_{n-1} (\varphi ) = c(\varphi (z',0) )$ and $(dd^c\varphi )^n = \delta_{ \{0\} }$.}
We set
$$\varphi_n (w_n) = \int\limits_{ \Delta^{n-1} } \varphi (z', w_n) (dd^c \varphi (z',w_n) )^{n-1}.$$
By Theorem 3.1 in \cite{ACKHZ09}, we have $\varphi_n\in\mathcal F(\Delta )$ and $dd^c\varphi_n = \delta_{ \{0\} }$. This implies that $\varphi_n (w_n) = \log |w_n|$. Take $c < c_{n-1} (\varphi ) = c(\varphi (z',0) )$ and $0<r<\frac 1 2 $ small enough such that
$$\int\limits_{\Delta_{2r}^{n-1} } e^{ -2c\varphi (z',0) } dV_{2n-2}(z') < +\infty.$$
By the effective version of the semicontinuity theorem for log canonical thresholds in \cite{DK01} (see also \cite{Hie14}), there exist $B,\delta >0$ such that
$$\int\limits_{\Delta_r^{n-1} } e^{ -2c\varphi (z',w_n) } dV_{2n-2} (z') \leq B,\ \forall |w_n |<\delta.$$
Using Lemma 2.4 for $z'\in\Delta^{n-1}\to \varphi (z',w_n)\in\mathcal E_1 (\Delta^{n-1} )$, we get
$$\aligned
& \int\limits_{\Delta_r^{n-1} } e^{ -2\lambda \varphi (z',w_n)} dV_{2n-2} (z')\\ 
&\leq \pi^{n-1} r^{2(n-1)} + \frac { B } { 2(\lambda-c) } e^{ 2(\lambda - c) |\varphi_n (w_n)|  c^{n-1} (n-1)^{-n+1}  } \\
&\ + c_{n-1}  \left(1 + 2^{n} |\varphi_n (w_n) | \lambda^{n} (n-1)^{-n} \right)^{n - 2}\frac { 1 } { 2(n-1) } e^{ 2(\lambda - c) |\varphi_n (w_n) | c^{n-1} n^{-n+1}  }\\
&\ + c_{n-1} \frac { (n-1)^2 } { \lambda (n+1) } \int\limits_{ 0 }^{ +\infty }  \left( 1 + 2x  \right)^{n - 2} e^ { -  2 (n-1)  x } dx.
\endaligned$$
This implies that
$$\lim\limits_{ w_n\to 0 }\int\limits_{\Delta_r^{n-1} } e^{-2\lambda\varphi (z',w_n )} dV_{n-2} (z') |w_n|^2 = 0,\ \forall \lambda\in (c,c+\frac { (n-1)^{n-1} } { c^{n-1} }).$$
By Lemma 2.2, we infer
$$c(\varphi ) \geq c+\frac { (n-1)^{n-1} } { c^{n-1} }.$$
Letting $c\to c_{n-1} (\varphi )$, we obtain
$$c(\varphi ) \geq c_{n-1} (\varphi ) + \frac { (n-1)^{n-1} } { c_{n-1}(\varphi )^{n-1} }.$$  
\qed 

\vskip1cm

{\bf Acknowledgment.} The author is grateful to Professor Jean-Pierre Demailly and Dr. Nguyen Ngoc Cuong for valuable comments. This paper was partly written when the author was visiting the Vietnam Institute for Advanced Study in Mathematics (VIASM). He would like to thank the VIASM for support and providing a fruitful research environment and hospitality. This research is funded by Vietnam National Foundation for Science and Technology Development (NAFOSTED) under grant number 101.02-2014.01.

\end{document}